\documentclass[12pt]{article}%
\usepackage[applemac]{inputenc}
\usepackage{amsmath,amssymb,fullpage}
\usepackage{amsthm}
\usepackage{amsfonts}
\usepackage{amsmath,amssymb}
\usepackage[applemac]{inputenc}
\usepackage{amsmath,amssymb,fullpage}
\usepackage{color}
\usepackage{color, colortbl}
\usepackage{amsmath}
\usepackage{amssymb}
\usepackage{graphicx}
\usepackage{tikz}
\usepackage{enumerate}
\usepackage{geometry}
\usetikzlibrary{arrows}
\usepackage{amsfonts}
\usepackage{amsmath,amssymb}
\usepackage[applemac]{inputenc}
\usepackage{amsmath,amssymb,fullpage}
\usepackage{color}
\usepackage{amsmath}
\usepackage{amssymb}
\usepackage{graphicx}
\usepackage{tikz}
\usepackage{graphicx}
\graphicspath{ {./images/} }
\usepackage{amsmath}
\usepackage{algpseudocode}
\usepackage{algorithm}
\usepackage{float}
\usepackage{mathabx}
\usepackage{siunitx}
\usepackage{placeins}

\let\Oldsection\section
\renewcommand{\section}{\FloatBarrier\Oldsection}

\let\Oldsubsection\subsection
\renewcommand{\subsection}{\FloatBarrier\Oldsubsection}

\let\Oldsubsubsection\subsubsection
\renewcommand{\subsubsection}{\FloatBarrier\Oldsubsubsection}

\usepackage{multirow}
\usepackage{caption}
\usepackage{subcaption}
\usepackage{bbm, dsfont}

\algdef{SE}[SUBALG]{Indent}{EndIndent}{}{\algorithmicend\ }%
\algtext*{Indent}
\algtext*{EndIndent}

\newtheorem{theorem}{Theorem}[section]
\newtheorem{lemma}[theorem]{Lemma}
\newcommand{\E}{E}

\newtheorem{corollary}[theorem]{Corollary}
\newtheorem{definition}[theorem]{Definition}

\newtheorem{remark}[theorem]{Remark}

\usepackage{tikz}
\usetikzlibrary{shapes.geometric, arrows}

\tikzstyle{input} = [circle, minimum width=1cm, text centered, draw=black, fill=green!20]

\tikzstyle{output} = [circle, minimum width=1cm, text centered, draw=black, fill=blue!20]

\tikzstyle{lstm} = [rectangle, rounded corners, minimum width=2cm, minimum height=1cm,text centered, draw=black, fill=red!20]

 \tikzstyle{lin} = [rectangle, minimum width=2cm, minimum height=1cm,text centered, draw=black, fill=orange!20]

 \tikzstyle{act} = [ellipse, minimum width=2cm, minimum height=1cm,text centered, draw=black, fill=yellow!20]

\tikzstyle{dot} = [rectangle, minimum width=2cm, minimum height=1cm,text centered]

\tikzstyle{arrow} = [thick,->,>=stealth]
\tikzstyle{map} = [thick, dashed,->,>=stealth]

\usetikzlibrary{arrows,shapes,automata,petri}

\newcommand{\dproof}{\noindent {Proof.} \quad}
\newcommand{\fproof}{\hfill $\square$ \bigskip}

\numberwithin{equation}{section}

\def\RB{\mathbb{R}}

\def\AC{\mathcal{A}}

\definecolor{LightCyan}{rgb}{0.88,1,1}

\def\RR{{\mathbb{ R}}}

\def\EE{{\mathbb{ E}}}

\def\1B{\text{1\!\!I}}

\def\R{\mathbb{R}}
\def\P{\mathbb{P}}
\def\E{\mathbb{E}}

\begin{document}

\title{Fokker-Planck equations for McKean-Vlasov SDEs driven by fractional Brownian motion}
\author{Saloua Labed $^{1}$, Nacira Agram $^{2}$ and Bernt Oksendal $^{3}$}
\date{9 January 2026}
\maketitle

\footnotetext[1]{Department of Mathematics, University of Biskra, Algeria.
\newline Email: s.labed@univ-biskra.dz}

\footnotetext[2]{Department of Mathematics, KTH Royal Institute of Technology
100 44, Stockholm, Sweden. \newline Email: nacira@kth.se. Work supported by
the Swedish Research Council grant (2020-04697).}

\footnotetext[3]{Department of Mathematics, University of Oslo, Norway.
\newline Email: oksendal@math.uio.no
}

\begin{abstract}
This paper investigates the probability distribution of solutions to McKean--Vlasov
stochastic differential equations driven by fractional Brownian motion with Hurst
parameter $H>\tfrac12$.

Our main contribution is the derivation of the associated Fokker--Planck equation,
which characterizes the time evolution of the law of the solution in a suitable
distributional framework. Under mild assumptions, we show that the law-valued
process is absolutely continuous in time and provide an explicit weak formulation
of the corresponding fractional McKean--Vlasov Fokker--Planck equation. In the case
where the law admits a density, we obtain a more explicit partial differential
equation with time-dependent diffusion coefficients induced by the fractional
noise.

We further establish a fractional Feynman--Kac representation, linking the forward
Fokker--Planck equation with a backward Kolmogorov equation for functionals of the
solution process. This result extends the classical Feynman--Kac framework to
mean--field dynamics driven by fractional Brownian motion.

To illustrate the theory, we analyze several explicit examples, including the law
of fractional Brownian motion itself and linear McKean--Vlasov fractional SDEs.
These examples highlight how fractional noise and mean--field interactions jointly
affect the probabilistic and analytic structure of the system.
\end{abstract}


\section{Introduction}
In this paper we study McKean-Vlasov stochastic differential equations (SDEs) driven
by fractional Brownian motion (fBm) $B^{(H)}$. The equations have the form
\begin{equation}
\begin{array}
[c]{l}%
dX\left(  t\right)  =\alpha\left(  t,X\left(  t\right)  ,\mu_{t}\right)
dt+\beta\left(  t,\mu_{t}\right)  {dB}^{(H)}\left(  t\right);  \quad 
X\left(  0\right)  =Z,
\end{array}
 \label{4.1}%
\end{equation}
where $\mu_{t}=\mathbb{P}(X(t))=\mathcal{L}(X(t))$ is the law of $X(t)$, and $B^{(H)}(t)$ is a
fractional Brownian motion with Hurst parameter $H\in(\tfrac{1}{2},1)$ on a
filtered probability space $(\Omega,\mathcal{F},\mathbb{F}=\{\mathcal{F}%
_{t}\}_{t\geq0},\mathbb{P})$. We assume that $Z$ is a given random variable in
$L^{2}(\mathbb{P})$, independent of $%
{{B}^{(H)}}$.\\

We explore the probability distributions of solutions of such McKean-Vlasov stochastic differential equations (SDEs) driven by fBm, extending the study of mean-field models under non-Markovian dynamics. McKean-Vlasov SDEs, which describe systems in which individual components evolve according to both their state and the distribution of the entire system, were rigorously studied by McKean \cite{Mc}. These equations have become central in various domains, including statistical physics, neuroscience, and economics, due to their ability to model large systems with interacting components by Carmona and Delarue \cite{CD}.

Fractional Brownian motion, introduced by Mandelbrot and Van Ness \cite{MN}, offers a generalization of classical Brownian motion, characterized by long-range dependence and self-similarity, making it suitable for modeling processes with memory and temporal correlations. The combination of fractional Brownian motion with McKean-Vlasov SDEs poses significant analytical challenges, as the non-Markovian nature of fBm complicates both the existence and uniqueness of solutions and the analysis of the corresponding probability distributions, see Fan et al. \cite{F} and for non McKean-Vlasov case, we refer to Hu \cite{Hu} and Biagini et al. \cite{BO, BHOZ}.

Our study begins by recalling key results by Fan et al. \cite{F}, on the existence and uniqueness of solutions for McKean-Vlasov SDEs driven by fBm. We then proceed to our main result, which is a proof of the associated Fokker-Planck equation, which describes the evolution of the probability density function of the system's solutions. In the case where the distribution is absolutely continuous with respect to the Lebesgue measure, we provide an explicit form of this Fokker-Planck equation.

To illustrate the applicability of our theoretical results, we solve several examples. These include the law of fractional Brownian motion and the geometric McKean-Vlasov SDE, which highlight the complex behavior induced by the combination of fractional noise and mean-field interactions. 

The paper is organized as follows.

In Section~2 we collect the analytic and stochastic preliminaries used throughout
the paper. This includes the introduction of a Hilbert space of measures, elements
of distribution theory, and the construction of fractional Brownian motion and
fractional white noise calculus within the white--noise framework.

In Section~3 we derive the Fokker--Planck equation associated with McKean--Vlasov
stochastic differential equations driven by fractional Brownian motion. We show that
the law of the solution satisfies a distributional evolution equation and obtain a
weak formulation for the case where the law admits a density, without imposing
classical $C^{1,2}$ regularity assumptions.

In Section~4 we establish a fractional Feynman--Kac representation for McKean--Vlasov
diffusions driven by fractional Brownian motion. This result provides the backward
Kolmogorov equation dual to the forward Fokker--Planck equation and connects the
probabilistic and PDE viewpoints.

In Section~5 we present explicit examples illustrating the theoretical results.
These include the recovery of the marginal law of fractional Brownian motion via a
time dependent heat equation, and the analysis of a linear geometric McKean--Vlasov
fractional SDE, for which both the Fokker--Planck equation and the Feynman--Kac
representation admit closed-form expressions.

\section{Preliminaries}

In this section we collect the analytic and stochastic preliminaries
needed throughout the paper. We first introduce a Hilbert space of
measures suitable for weak formulations of McKean--Vlasov equations.
We then recall elements of fractional Brownian motion and fractional
white noise calculus based on the white--noise framework.
All results presented here are standard and included to fix notation
and ensure the paper is self contained.

\subsection{A Hilbert space of measures}

Let $\mathbb{M}$ denote the space of (possibly random) finite measures
$\mu$ on $\mathbb{R}^d$ equipped with the norm
\begin{equation}\label{norm}
\|\mu\|_{\mathbb{M}}^{2}
:= \mathbb{E}\!\left[\int_{\mathbb{R}^d}
|\widehat{\mu}(y)|^{2}e^{-y^{2}}\,dy\right],
\end{equation}
where $\widehat{\mu}$ denotes the Fourier transform
\begin{equation}
\widehat{\mu}(y)
:= \int_{\mathbb{R}^d} e^{-i x\cdot y}\,\mu(dx),
\qquad y\in\mathbb{R}^d.
\end{equation}

For $\mu,\eta\in\mathbb{M}$ we define the inner product
\[
\langle \mu,\eta\rangle_{\mathbb{M}}
:= \mathbb{E}\!\left[
\int_{\mathbb{R}^d}
\operatorname{Re}\big(\overline{\widehat{\mu}(y)}\,\widehat{\eta}(y)\big)
e^{-y^{2}}\,dy\right].
\]
The completion of this space is again denoted by $\mathbb{M}$, and
$\mathbb{M}_0$ denotes its deterministic elements.

In particular, if $\mu(dx)$ is absolutely continuous with respect to Lebesgue measure $dx$ with Radon-Nikodym-derivative $m(x)=\frac{\mu(dx)}{dx}$, so that $\mu(dx)=m(x) dx$ with $m \in L^1(\RR^d)$, we define
the Fourier transform
of $m$ at $y$, denoted by $\widehat{m}(y)$ or $F[m](y)$, by
\begin{align}
F[m](y)=\widehat{m}(y)=\int_{\mathbb{R}^d} e^{-ixy} m(x) dx; \quad y \in \mathbb{R}^d.
\end{align}

\subsection{Schwartz space of tempered distributions}
Let $\mathcal{S}=\mathcal{S}(\mathbb{R})$\label{simb-028} be \textit{the
Schwartz space of rapidly decreasing smooth $C^{\infty}(\mathbb{R}^d)$ real
functions} 
on $\mathbb{R}$. The space $\mathcal{S}=\mathcal{S}(\mathbb{R})$ is a
Fr\'echet space with respect to the family of semi-norms:\label{simb-029} 
\begin{equation*}
\Vert f \Vert_{k,\alpha} := \sup_{x \in \mathbb{R}}\big\{ (1+|x|^k) \vert
\partial^{{n}}  f(x)\vert \big\},
\end{equation*}
where $k,{{n}} = 0,1,...$

Let $\mathcal{S}^{\prime }=\mathcal{S}^{\prime }(\mathbb{R})$\label%
{simb-031} be its dual, called the space of \emph{tempered distributions}. 
\index{tempered distributions} Let $\mathcal{B}$ denote the family of all
Borel subsets of $\mathcal{S}^{\prime }(\mathbb{R})$ equipped with the
weak {*} topology. If $\Phi \in \mathcal{S}^{\prime }$ and $f \in \mathcal{%
S}$ we let \label{simb-033} 
\begin{equation}
\Phi (f) \text{ or } \langle \Phi ,f \rangle  \label{3.1}
\end{equation}%
denote the action of $\Phi $ on $f $. For example, if $\Phi =m$ is a
measure on $\mathbb{R}$ then 
\begin{equation*}
\langle \Phi,f \rangle =\int\limits_{\mathbb{R}^{d}} f(x)dm(x),
\end{equation*}%
and, in particular, if this measure $m$ is concentrated on $x_{0}\in \mathbb{%
R}$, then 
\begin{equation*}
\langle \Phi , f \rangle = f(x_{0})
\end{equation*}
is evaluation of $f $ at $x_{0}\in \mathbb{R}$.\newline
Other examples include 
\begin{equation*}
\left\langle \Phi , f \right\rangle = f^{\prime }(x_{1}),
\end{equation*}%
i.e. $\Phi $ takes the derivative of $f$ at a point $x_{1}$.%
\newline
Or, more generally, 
\begin{equation*}
\left\langle \Phi , f \right\rangle = f^{(k)}(x_{k}),
\end{equation*}%
i.e. $\Phi $ takes the $k$'th derivative at the point $x_{k}$,
or linear combinations of the above.\newline
If $\Phi \in \mathcal{S}'$ we define its Fourier transform $\widehat{\Phi} \in \mathcal{S}'$
by the identity
\begin{align}
\langle \widehat{\Phi}, f \rangle = \langle \Phi, \widehat{f} \rangle; \quad f \in \mathcal{S}.
\end{align}
The partial derivative with respect to $x_k$ of a tempered distribution $\Phi$ is defined by 
\begin{align*}
&\langle \tfrac{\partial}{\partial x_k} \Phi, f \rangle = -\langle \Phi, \tfrac{\partial}{\partial x_k} f \rangle; \quad \phi \in \mathcal{S}^d.\\
&\text{ More generally, }\\
&\langle \partial^{\alpha} \Phi, f \rangle = (-1)^{|\alpha|} \langle \Phi, \partial^{\alpha} f \rangle; \quad \phi \in \mathcal{S}^d.
\end{align*}
We refer to Chapter 8 in {Folland} \cite{F} for more information.\\

\subsection{fBm and its calculus}
In this subsection we briefly review the basic notation and results of the white noise probability space and its associated Brownian motion. 
Let $\mathcal{S}(\RR)$ denote the Schwartz space of rapidly decreasing smooth functions on $\mathbb{R}$ and define $\Omega={\cal S}'(\RR)$ to be its dual, which is the Schwartz space of tempered distributions, equipped
with the weak-star topology. This space will be the base of our basic \emph{white noise
probability space}, which we recall in the following:
\vskip 0.3cm
As events we will use the family ${\cal
B}({\cal S}'(\RR))$ of Borel subsets of ${\cal S}'(\RR)$, and our
probability measure {$\P$} is defined as follows:

\begin{theorem}{\bf (The Bochner--Minlos theorem)}\\
There exists a unique probability measure $\P$ on ${\cal B}({\cal S}'(\R))$
with the following property:
$$\E[e^{i\langle\cdot,\phi\rangle}]:=\int\limits_{\cal S'}e^{i\langle\omega,
\phi\rangle}d\mu(\omega)=e^{-\tfrac{1}{2} \Vert\phi\Vert^2};\quad i=\sqrt{-1}$$
for all $\phi\in{\cal S}(\R)$, where
$\Vert\phi\Vert^2=\Vert\phi\Vert^2_{L^2(\RR^n)},\quad\langle\omega,\phi\rangle=
\omega(\phi)$ is the action of $\omega\in{\cal S}'(\R)$ on
$\phi\in{\cal S}(\R)$ and $\E=\E_{\P}$ denotes the expectation
with respect to  $\P$.
\end{theorem} 

We will call the triplet $({\cal S}'(\R),{\cal
B}({\cal S}'(\R)),\P)$ the {\it  white noise probability
space\/}, and $\P$ is called the {\it white noise probability measure}.

The measure $\P$ is also often called the (normalised) {\it Gaussian
measure\/} on ${\cal S}'(\R)$. It is not difficult to prove that if $\phi\in L^2(\R)$ and
we choose $\phi_k\in{\cal S}(\R)$ such that $\phi_k\to\phi$ in $L^2(\R)$,
then
$$\langle\omega,\phi\rangle:=\lim\limits_{k\to\infty}\langle\omega,\phi_k\rangle
\quad\text{exists in}\quad L^2(\P)$$
and is independent of the choice of $\{\phi_k\}$. In
particular, if we define
$$\widetilde{B}(t):=\widetilde{B}(t,\omega)=\langle\omega,\chi_
{[0,t]}\rangle; $$ 
then $\widetilde{B}(t,\omega)$ has a $t$-continuous version $B(t,\omega)$, which
becomes a { \emph{Brownian motion}}, in the following sense:
 
By a  \emph{Brownian motion} we mean a family
$\{X(t,\cdot)\}_{t\in\R_{+}}$ of random variables on a probability space
$(\Omega,{\cal F},\P)$ such that
\begin{itemize} 
\item
$X(0,\cdot)=0\quad\text{almost surely with respect to } \P,$
\item
$\{X(t,\omega)\}_{t\geq 0,\omega \in \Omega}$ is a continuous and Gaussian stochastic process 
\item
For all $s,t \geq 0$,
$\E[X(s,\cdot) X(t,\cdot)]=\min(s,t).$ 
\end{itemize} 

It can be proved that the process $\widetilde{B}(t,\omega)$ defined above has a modification $B(t,\omega)$ which satisfies all these properties.
This process $B(t,\omega)$ then becomes a Brownian motion. This is the version of Brownian motion we will use to construct the associated fBm, as explained in the next subsection.

\subsection{The operator $M$}
\begin{definition} \label{fBm}
By a \emph{fBm} with Hurst parameter $H \in (0,1)$ we mean a continuous Gaussian stochastic process $\{Y(t)\}_{t\geq 0}$ with mean 0 and covariance
$$\E[Y(s) Y(t)]=C(|s|^{2H} + |t|^{2H}-|t-s|^{2H} ); \text{ for all } s,t\geq 0,$$
for some constant $C>0$.
\end{definition}

fBm can be constructed in several ways. We will follow the approach used by, among others,  Elliot and van der Hoek \cite{EvdH}, where fBm is constructed by applying an operator $M$ to the classical Bm. See also \cite{L}, \cite{ST} and \cite{Be1},\cite{Be2}. This approach has the advantage that all the different fBm's corresponding to different Hurst coefficients $H \in (0,1)$ are defined on the same white noise probability space $(\Omega, \mathbb{F},P)$ introduced below. Moreover, it has several technical advantages, as we have exploited in the proof of our main result (Theorem 3.1).\\

\begin{definition}(\cite{EvdH},(A10))
Let $F^{-1}$ denote the inverse Fourier transform, defined by
 \begin{align}
 F^{-1} g(x):= \frac{1}{\sqrt{2 \pi}}\int_{\mathbb{R}} e^{ixy} g(y) dy,
 \end{align}
 for all functions $g$ such that the integral converges.

Let $H \in (0,1)$ and let $\mathcal{S}$ be the Schwartz space of rapidly decreasing smooth functions on $\mathbb{R}$. The operator $M=M^{(H)}: \mathcal{S} \mapsto \mathcal{S}$ is defined by
\begin{align}
Mf(x)&=F^{-1}(|y|^{\frac{1}{2}-H} (Ff)(y))(x)\nonumber\\
&=
\textstyle{\frac{1}{\sqrt{2 \pi}}\int_{\mathbb{R}} e^{ixy} (|y|^{\frac{1}{2} - H}\frac{1}{\sqrt{2 \pi}} \int_{\mathbb{R}} e^{-izy}f(z)dz) dy}
\end{align}
If $H > \frac{1}{2} $ then the operator $M$ can also be defined as follows:
\begin{align}
M_Hf(x)= C_H \int_{\mathbb{R}} \frac{f(y+x)}{| y | ^{\frac{3}{2} - H}} dy; \quad f \in \mathcal{S},\label{MH}
\end{align}
where
\begin{equation}
C_H=[2\Gamma(H-\frac{1}{2}) \cos(\frac{\pi}{2}(H-\frac{1}{2}))] ^{-1},
\end{equation}
and $\Gamma$ is the classical Gamma function.
\end{definition}
\emph{From now on we assume that $H \in (\frac{1}{2},1)$.}\\

Since now $H \in (\frac{1}{2},1)$ is fixed, we will suppress the index $H$ and write $M^{(H)} = M$ in the following.\\
We can in a natural way extend the operator $M$  from $\mathcal{S}$ to the space
\begin{align}
L_{H}^2(\mathbb{R}):= \{ f: \mathbb{R} \mapsto \mathbb{R}; Mf \in L^2(\mathbb{R} )\}.
\end{align}
Note that $L_{H}^2 (\mathbb{R})$ is a Hilbert space when equipped with the inner product
\begin{equation}
(f,g)_H := (f,g)_{L_H^2(\mathbb{R})}=(Mf,Mg)_{L^2 (\mathbb{R})};\quad f,g \in L^2_{H}(\mathbb{R})
\end{equation}
and the corresponding norm
\begin{equation}
|f|_H=\Big((Mf,Mf)_{L^2 (\mathbb{R})}\Big)^{\tfrac{1}{2}};\quad f \in L^2_{H}(\mathbb{R}).
\end{equation}

We will need the following results about the operator $M$:

\begin{lemma}(\cite{DO} )\label{M2} 
The operator $M$ satisfies following property:
\begin{align} \label{eq3.5}
M^2f(x)= C_H^2 \int_{\mathbb{R}} \int_{\mathbb{R}} |yz|^{H-\frac{3}{2}} f(y+z+x)dy dz; \quad f \geq 0.
\end{align}
\end{lemma}
\begin{lemma} \cite{DO}
Let  $\varphi,\psi:\mathbb{R}\rightarrow \mathbb{R}^{+}$. Then
	\begin{align} \label{eq l21}
	\int_{\mathbb{R}}\varphi(t)	M^2\psi(t)dt= \int_{\mathbb{R}}\psi(t)	M^2\varphi(t)dt=\int_{\mathbb{R}} M\varphi(t) M\psi(t) dt.
	\end{align}
\end{lemma}

\subsection{Fractional Brownian motion (fBm)}

We can now give our definition of fBm:

\begin{definition}{(Fractional Brownian motion)}\\
For $t \in \mathbb{R}$ define
\begin{equation} \label{eq4.5a}
B_t^{(H)}:= B^{(H)}(t):= B^{(H)}(t,\omega):=\langle \omega, M_t(\cdot) \rangle ; \quad \omega \in \Omega 
\end{equation}
where we for simplicity of notation put 
\begin{align}
M_t (x):=M(\chi_{[0,t]})(x)=C_H \Big(\frac{t-x}{|t-x|^{\frac{3}{2}-H}}+\frac{x}{|x|^{\frac{3}{2}-H}}\Big);\quad t,x \in \mathbb{R}.
\end{align}
See (A7), p. 324 in \cite{EvdH}.
Then $B^{(H)}$ is a Gaussian process, and  we have 
\begin{equation}
\E[B^{(H)} (t)]=B^{(H)}(0)=0,
\end{equation}
and 
\begin{equation}
\E[B^{(H)}(s) B^{(H)}(t)]= C_H[ |t|^{2H}+ |s|^{2H}-|t-s|^{2H}]; \quad s,t \in \mathbb{R}. \label{var}
\end{equation}
(See \cite{EvdH},(A10).)

Hence 
\begin{align}
\E[B^{(H)}(t)^2]= 2 C_H t^{2H}=:\sigma^2, \label{var}
\end{align}\label{2.19}
and from this it follows, since $B^{(H)}(t)$ is Gaussian,
\begin{align}
\E[(B^{(H)}(t))^4]= 3 \sigma^4 = 12 C_H^2 t^{4H}.
\end{align} \label{moment}


It follows from the Kolmogorov continuity theorem that  the process $B^{(H)}$ has a continuous version, which we will also denote by $B^{(H)}.$  This process satisfies all the conditions of being an fBm given in Definition \ref{fBm} and we will use this as our
 \emph{ the fBm with Hurst parameter $H$}.
\end{definition}

\begin{remark}{\bf{(A remark about filtrations)}}
We let $\mathcal{F}=\{\mathcal{F}_t\}_{t\geq 0}$ denote the filtration generated by Brownian motion $B(t,\omega)$. This means that, for all $t$,  $\mathcal{F}_t$ is the $\sigma$-algebra generated by the random variables $B(s,\cdot); s \leq t$. Similarly $\mathcal{F}^{(H)}=\{\mathcal{F}_t^{(H)}\}_{t\geq 0}$ denotes the filtration generated by $B^{(H)}(t,\omega)$.
Recalling that $B(t,\omega)=\langle \omega, \chi_{[0,t]} \rangle, $ while $B^{(H)}(t,\omega)=\langle \omega,M_t \rangle$, we see that the two filtrations are not identical. In particular, $B^{(H)}(t)$ is not $\mathbb{F}$-adapted.\\
This is in contrast with the situation for some other constructions of fBm, where the filtration of $B(t)$ and $B^{(H)}(t)$ are the same. See e.g. \cite{BHOZ}.
\end{remark}

\subsection{Fractional white noise calculus}
Using that
\begin{equation}
B^{(H)}(t)=\int_{\mathbb{R}} \chi_{[0,t]}(s) dB^{(H)}(s)= \int_{\mathbb{R}}M_t(s)dB(s)=\langle \omega, M_t \rangle.\label{fBm2}
\end{equation}
We see that if $f$ is a simple function of the form $f=\sum c_k \chi_{[t_k, t_{k+1}]}$ then
\begin{align}
\int f(t) dB^{(H)}(t)&= \sum c_k (B^{(H)}(t_{k+1} )- B^{(H)}(t_k))=\sum c_k \langle \omega, M_{t_{k+1}} - M_{t_k} \rangle \nonumber\\
&=\langle \omega, \sum c_k (M_{t_{k+1} }-M_{t_k} \rangle 
=\langle \omega,M(\sum c_k \chi_{[t_k,t_{k+1}]}\rangle \nonumber\\
&= \langle \omega,Mf \rangle,
\end{align}
where $\chi$ denotes the indicator function.
From this it follows by approximating $f \in \mathcal{S}$ with step functions that
\begin{equation}
 \int_{\mathbb{R}} f(t) dB^{H}(t)=\langle \omega,Mf \rangle =\int_{\mathbb{R}} Mf(t) dB(t),
\end{equation}
for all $f \in \mathcal{S}$.
\begin{remark}
    Note that \eqref{fBm2} combined with Theorem 8.5.2 in \cite{O} implies that there is a Brownian motion  $\widetilde{B}(t)$ such that
    \begin{align}
        B^{(H)}(t)= \widetilde{B}(\rho(t)); \text{ where } \rho(t) = \int_0^t M_t(s)^2 ds.
    \end{align}
    In other words, $B^{(H)}(t)$ has the same law as a time-changed Brownian motion, with time change $\rho(t) = \int_0^t M_t(s)^2 ds.$
\end{remark}
We know that the action of $\omega \in \mathcal{S}'(\mathbb{R})$ extends from $\mathcal{S}(\mathbb{R})$ to $L^2(\mathbb{R}),$ and that if $f \in \mathcal{S}$ then $Mf \in L^2(\mathbb{R}).$ Therefore we can define $M$ on $\mathcal{S}'(\mathbb{R})$ by setting
\begin{equation}
\langle M\omega ,\varphi \rangle= \langle \omega, M\varphi \rangle; \varphi \in \mathcal{S}(\mathbb{R}), \omega \in \mathcal{S}'(\mathbb{R}).
\end{equation}
Using this, we can write

\begin{align}
B^{(H)}(t,\omega)=\langle \omega, M_t \rangle= \langle M\omega, \chi_{[0.t]} \rangle =B(t,M\omega).
\end{align}

We define the \emph{fractional white noise} 
$W^{H}(t)=W^{H}(t,\omega)$ by
\begin{align}
W^{H}(t,\omega)=\frac{d B^{(H)}(t,\omega)}{dt} \text{ (derivative in } (\mathcal{S})^{*} ).
\end{align}
By the above it follows that 
\begin{align}
W^{H}(t,\omega)=\frac{d B(t,M\omega)}{dt}=W(t,M\omega),
\end{align}
 where 
\begin{align}
W(t)=W(t,\omega)=\frac{d B(t,\omega)}{dt} \text{ is the classical white noise in } (\mathcal{S})^{*}.
\end{align}
For $0 \leq T \leq \infty$ let $\mathcal{L}_H^2[0,T]$ denote the set of measurable processes $\varphi(t)=\varphi(t,\omega)$ such that 
\begin{align}
\E[\int_0^T (M \varphi)^2(t)dt] < \infty.
\end{align}
Then we define the \emph{Wick-Ito-Skorohod} (WIS)-integral  of $\varphi \in \mathcal{L}_H^2[0,T]$  as follows:

\begin{definition}  \label{wis}
(The  Wick--It\^o--Skorohod (WIS) integral) \\
Let  $\varphi : \mathbb{R} \rightarrow (\mathcal{S})^{*}$ be such that $\varphi (s)  \diamond W^{H}(s) $ is $ds$-integrable in $ (\mathcal{S})^{*}$. Then we say that  $\varphi$ is
$dB^{(H)}$-integrable and we define the Wick--It\^o--Skorohod (WIS) integral of $\varphi (\cdot) = \varphi (\cdot, \omega)$ with respect to $B^{(H)}$ by
\begin{align}
\int_{\mathbb{R}} \varphi(s) dB^{H}(s)= \int_{\mathbb{R}} \varphi(s) \diamond W^{H}(s) ds.
\end{align} 
\end{definition}

Let $\mathcal{L}_H^2[0,T]$ denote the set of measurable processes $\varphi$ such that 
\begin{align}
    \E[\int_0^T (M\varphi)^2(t)dt] < \infty.
\end{align}We recall the following theorem from \cite{EvdH}:
\begin{theorem}
Let $\varphi \in \mathcal{L}_H^2[0,T].$ \\
Then
\begin{equation} 
\int_{\mathbb{R}} \varphi(s) dB^{(H)}(s)=\int_{\mathbb{R}} M\varphi(s) dB(s)=\int_{\mathbb{R}} M\varphi(s) \diamond W(s) ds,
\end{equation} 
where $\diamond$ denotes the Wick product, the second integral is the classical Wick-Ito-Skorohod (WIS)-integral for $B$ and the last integral is the (Pettis-) Lebesgue integral in the Hida space $(\mathcal{S})^{*}$ of stochastic distributions.
\end{theorem}

The next lemma states that the expectation of a WIS-integral with respect to a fractional Brownian motion is zero:

\begin{lemma}  Let $\varphi \in \mathcal{L}_H^2[0,T]$. Then
	$$\E\left[\int_0^{T}  \varphi(s) dB^{(H)}(s) \right]=0.$$
\end{lemma}
\dproof This follows from Definition \ref{wis} and the fact that  the expectation of a Wick product is the product of the expectations.
\fproof

We also note the following relation between the Wick product $\diamond$ and the ordinary product:
\begin{lemma}\cite{EvdH}
Let $ f, g \in \mathcal{L}_H^2[0,T].$ Then
\begin{align} \label{iso}
&(\int_0^t f(s) dB^{(H)}(s)) \diamond (\int_0^t g(s) dB^{(H)}(s))\nonumber\\
&=(\int_0^t f(s) dB^{(H)}(s)) (\int_0^t g(s) dB^{(H)}(s)) -(f,g)_{H}.
\end{align}
\end{lemma}

We need the following It\^o formula  for fractional WIS processes: 
\begin{theorem}\cite{BO} \label{Ito}{\bf{(An It\^o formula for fractional WIS processes)}}\newline
Let $H > \frac{1}{2}$. Let $\alpha$ be a measurable process such that $\E[\int_0^t |\alpha(s)|ds] <\infty$ for all $t \geq 0$ and let $\beta$ be an $\mathbb{F}$-adapted cadlag WIS-integrable process. Suppose $X(t)$ is a \emph{fractional WIS process} of the form
\begin{align}
dX(t)=\alpha(t)dt + \beta(t) dB^{(H)}(t);\quad X(0)=x \in \mathbb{R},
\end{align}
Let $f:[0,T] \times\mathbb{R} \mapsto \mathbb{R}$ be in  $C^{1,2}$ and put $Y(t)=f(t,X(t))$. Then
\begin{align}
dY(t)&=\frac{\partial}{\partial t}f(t,X(t)) dt +\frac{\partial}{\partial x}f(t,X(t)) dX(t) \nonumber\\
&+ \frac{ \partial^2}{\partial x^2} f(t,X(t))\beta(t) M^2(\beta \chi_{[0,t]})(t)dt.
\end{align}
\end{theorem}

\dproof
Since $\beta$ is adapted we have $MD_t^{(H)} \beta(s)=D_t \beta(s)=0$ for all $s < t$. Therefore this result follows as a special case of Theorem 2.11 in \cite{BO}.
\fproof 

An important special case of the (2-dimensional) It\^o formula is the following  product rule:
\begin{lemma} \label{add2} (Product rule)
Let $\theta_i,\sigma_i; i=1,2$ be as $\theta$ and $\sigma$ respectively in Theorem \ref{Ito} and let
	
	$$dX_i=\theta_i(t)dt+\sigma_i(t)dB^{(H)}(t), \quad i=1,2$$
	with $X_i(0)=x_i.$ Then the following holds:
	\begin{eqnarray}
	&&\E[X_1(t)X_2(t)]=x_1x_2+\E\Big[\int_0^t\Big( X_1(s)\theta_2(s)+X_2(s)\theta_1(s)\Big)ds\Big]\label{8}\\
	&&\phantom{\E[X_1(t)X_2(t)]=}+\int_0^t\Big(\sigma_1(s)M^2(\sigma_2\chi_{[0,t]})(s)+\sigma_2(s)M^2(\sigma_1\chi_{[0,t]})(s)\Big) ds \Big]. \nonumber
	\end{eqnarray}
\end{lemma}

\subsection{Well-posedness of the McKean-Vlasov fBm equation}

Throughout the paper we consider McKean-Vlasov stochastic differential
equations driven by fBm of the form
\begin{equation}\label{eq:MK}
dX(t)
 = \alpha(t,X(t),\mu_t)\,dt
 + \beta(t,\mu_t)\,dB^{(H)}(t), 
 \qquad X(0)=Z,
\end{equation}
where $\mu_t=\mathcal{L}(X(t))$ and $B^{(H)}$ is an fBm with Hurst index
$H>1/2$.

We impose the standard Lipschitz assumptions:
\begin{description}
\item[(H)] There exists a nondecreasing function $K:[0,T]\to\R_+$ such that
for all $t\in[0,T]$, $x,y\in\R$ and $\mu,\nu\in \mathbb{M}_0$,
\[
\begin{aligned}
|\alpha(t,x,\mu)-\alpha(t,y,\nu)|
 &\le K(t)\big(|x-y|+\|\mu-\nu\|_{\mathbb{M}_0}\big),\\
|\beta(t,\mu)-\beta(t,\nu)|
 &\le K(t)\,\|\mu-\nu\|_{\mathbb{M}_0},
\end{aligned}
\]
and
\[
|\alpha(t,0,\delta_0)| + |\beta(t,\delta_0)| \le K(t).
\]
\end{description}

Under these assumptions, the McKean-Vlasov equation \eqref{eq:MK} is
well posed:

\begin{theorem}[Well-posedness; see \cite{DO}, Th.\,4.1]
If $Z\in L^2(\Omega)$ and {\rm(H)} holds, then the equation
\eqref{eq:MK} admits a unique solution 
$X\in\mathcal S^2([0,T])$, and the law map
$t\mapsto\mu_t=\mathcal{L}(X(t))$ is continuous in $\mathbb{M}_0$.
\end{theorem}

We emphasize that the existence and uniqueness theory is not the focus of
this paper; we rely on the well-posedness results of \cite{DO} and
concentrate instead on the analysis of the associated
Fokker-Planck equation and the fractional Feynman-Kac formula.

\begin{remark}
From \eqref{2.19} we have
\[
\E[(B^{(H)}(t))^2]=2C_H t^{2H}.
\]
On the other hand, applying the product rule
(Lemma~\ref{add2}) with $X_1=X_2=B^{(H)}$ yields
\[
\E[(B^{(H)}(t))^2]
=
2\int_0^t M^2(\chi_{[0,s]})(s)\,ds.
\]
Differentiating with respect to $t$ gives the following useful result:
\begin{align}
M^2(\chi_{[0,t]})(t)
=\frac12\frac{d}{dt}\E[(B^{(H)}(t))^2] \label{M2}
=2H\,C_H\,t^{2H-1}.
\end{align}

\end{remark}

\begin{lemma}[Absolute continuity of the law process for McKean--Vlasov fBm SDE]\label{lem2.14}
\label{lem:AC-fractional}
Let $X=(X_t)_{t\in[0,T]}$ be the solution of the McKean--Vlasov fBm equation
\begin{equation}\label{MK-again}
dX(t)
= \alpha(t,X(t),\mu_t)\,dt + \beta(t,\mu_t)\,dB^{(H)}(t), 
\qquad X(0)=Z,
\end{equation}
with $H> \frac12$ and $\mu_t=\mathcal{L}(X_t)$.
Assume that:
\begin{itemize}
  \item[(i)] $\alpha:[0,T]\times \R\times \mathbb{M}_0\to\R$ and 
             $\beta:[0,T]\times \mathbb{M}_0\to\R$ are bounded
             and continuous in $t$ (and in $x$ for $\alpha$) for each fixed measure;
  \item[(ii)] the fractional variance factor $M^2(\beta\chi_{[0,t]})(t)$ is bounded on $[0,T]$, i.e.,
  \[
    \sup_{t\in[0,T]}
    \big| M^2(\beta\chi_{[0,t]})(t) \big|
    < \infty;
  \]
  \item[(iii)] for each $t\in[0,T]$, the law $\mu_t$ belongs to the Hilbert space $\mathbb{M}_0$.
\end{itemize}
Then the map
\[
t \longmapsto \mu_t \in \mathbb{M}_0
\]
is absolutely continuous on $[0,T]$. In particular,
\[
\frac{d}{dt}\mu_t := \lim_{h\to 0}\frac{\mu_{t+h}-\mu_t}{h}
\]
exists in $\mathbb{M}_0$ for all $t\in[0,T]$, and we have the estimate
\begin{equation}\label{eq:AC-bound}
\|\mu_{t+h}-\mu_t\|_{\mathbb{M}_0} \le C\,|h|,
\qquad t,t+h\in[0,T],
\end{equation}
for some constant $C>0$ independent of $t$ and $h$.
\end{lemma}
\dproof
Fix $t\in[0,T]$ and $h$ such that $t,t+h\in[0,T]$.
Let $\mu_t=\mathcal{L(}X_t)$ and denote
\begin{align}
\widehat{\mu}_t(y) := F[\mu_t](y)
= \E\big[e^{-iyX_t}\big],\qquad y\in\R.
\end{align}
Define
\begin{align}
\varphi_y(x) := e^{-iyx}, \qquad x,y\in\R. \label{test}
\end{align}
By the fractional It\^o formula (Theorem~\ref{Ito}) applied to $\varphi_y(X_s)$ on $[t,t+h]$, we have
\[
\E\big[\varphi_y(X_{t+h}) - \varphi_y(X_t)\big]
= \E\Big[\int_t^{t+h} A_s\varphi_y(X_s)\,ds\Big],
\]
where the time-dependent generator $A_s$ acts on $\varphi\in C^2(\R)$ as
\[
A_s\varphi(x)
= \alpha(s,x,\mu_s)\,\varphi'(x)
 + \beta(s,\mu_s)\,M^2(\beta\chi_{[0,s]})(s)\,\varphi''(x).
\]
For our choice $\varphi_y(x)=e^{-iyx}$ we obtain
\[
A_s\varphi_y(x)
= \Big(-iy\,\alpha(s,x,\mu_s)
      -y^2\,\beta(s,\mu_s)\,M^2(\beta\chi_{[0,s]})(s)\Big)e^{-iyx}.
\]
Since $\alpha$ and $\beta$ are bounded and $M^2(\beta\chi_{[0,s]})(s)$ is bounded in $s$ by assumption, there exists a constant $C_1>0$ such that
\[
\big|A_s\varphi_y(x)\big|
\le C_1\,(|y| + y^2)\,\big|e^{-iyx}\big|
\le C_2\,(1+y^2),
\]
for all $s\in[0,T]$, $x\in\R$ and $y\in\R$, where we used $|e^{-iyx}|=1$ and absorbed constants into $C_2$.\\
Taking expectations and integrating in time, we get
\begin{align*}
\big|\widehat{\mu}_{t+h}(y) - \widehat{\mu}_t(y)\big|
&= \left|\E\big[\varphi_y(X_{t+h}) - \varphi_y(X_t)\big]\right| \\
&= \left|\int_t^{t+h} \E\big[A_s\varphi_y(X_s)\big]\,ds\right| \\
&\le \int_t^{t+h} \E\big[|A_s\varphi_y(X_s)|\big]\,ds
 \le |h|\,C_2\,(1+y^2).
\end{align*}
Therefore,
\begin{align*}
\|\mu_{t+h}-\mu_t\|_{\mathbb{M}_0}^2
&= \int_{\R} \big|\widehat{\mu}_{t+h}(y) - \widehat{\mu}_t(y)\big|^2 e^{-y^2}\,dy \\
&\le |h|^2 C_2^2 \int_{\R} (1+y^2)^2 e^{-y^2}\,dy
=: C_3^2\,h^2,
\end{align*}
where the integral is finite since $(1+y^2)^2 e^{-y^2}$ is integrable.\\
Taking square roots yields
\[
\|\mu_{t+h}-\mu_t\|_{\mathbb{M}_0}
\le C_3\,|h|,
\]
which is exactly \eqref{eq:AC-bound} with $C:=C_3$.\\
Thus $t\mapsto\mu_t$ is Lipschitz (hence absolutely continuous) as a map from $[0,T]$ into the Hilbert space $\mathbb{M}_0$. By standard results on absolutely continuous curves in Hilbert spaces, the strong derivative
\[
\mu'_t := \frac{d}{dt}\mu_t
= \lim_{h\to 0}\frac{\mu_{t+h}-\mu_t}{h}
\]
exists in $\mathbb{M}_0$ for all $t\in[0,T]$, which completes the proof.
\fproof

\begin{remark}
Assumption~(ii) in Lemma \ref{lem2.14} is satisfied in particular when $\beta$ is constant.
Indeed, \eqref{M2} yields
$M^2(\beta\chi_{[0,t]})(t)=\beta 2H C_Ht^{2H-1}<\infty$ for all $t\in[0,T]$.
\end{remark}

\section{A Fokker--Planck Equation for McKean--Vlasov SDEs Driven by fBm}
Let $X(t)=X_t\in\mathbb{R}$ solve the mean-field SDE driven by a fBm
(henceforth called a \emph{McKean--Vlasov fBm equation})
\begin{equation}\label{MK}
dX(t)=\alpha(t,X(t),\mu_t)\,dt+\beta(t,\mu_t)\,dB^{(H)}(t),
\qquad X(0)=Z,
\end{equation}
where $\mu_t=\mathcal{L}(X(t))$ denotes the law of $X(t)$ and $H>1/2$.
For simplicity we assume that
\begin{itemize}
\item $\alpha:[0,T]\times\mathbb{R}\times\mathbb{M}_0\to\mathbb{R}$  
      and 
      $\beta:[0,T]\times\mathbb{M}_0\to\mathbb{R}$  
      are bounded functions;
\item $\alpha$ is continuous in $(t,x)$ for each fixed $\mu$;
\item $\beta$ is continuous in $t$ for each fixed $\mu$.
\end{itemize}
In this section we derive the evolution equation for the law $\mu_t$.  
The resulting PDE is the fractional McKean--Vlasov Fokker--Planck equation.
\begin{theorem}[Fokker--Planck equation]\label{FP1}
Assume that $t\mapsto \mu_t$ is continuous in $\mathcal{S}'=\mathcal{S}'(\mathbb{R})$.
Then $(\mu_t)_{t\in[0,T]}$ is a solution of the Fokker--Planck equation
\begin{equation}\label{eq:FP-weak}
\frac{d}{dt}\mu_t = A^*\mu_t \quad \text{in } \mathcal{S}'.
\end{equation}

Equivalently,
\begin{equation}\label{eq:FP-integral}
\mu_t = \mu_0 + \int_0^t A^*\mu_s \, ds\quad\text{in }\mathcal{S}'.
\end{equation}
Here
\begin{align}
A^{*}\mu_t
&=-D\!\left[\alpha(t,\cdot,\mu_t)\,\mu_t(\cdot)\right]
+\beta(t,\mu_t)\,M^2(\beta\chi_{[0,t]})(t)D^2[\mu_t(\cdot)]\, \label{A*}\\
&\text{ which is the dual of }  \nonumber\\
A\mu_t&= \alpha (t,\cdot,\mu_t(\cdot) D[\mu_t(\cdot)] +\beta(t,\mu_t)\,M^2(\beta\chi_{[0,t]})(t)
D^2\!\left[\mu_t(\cdot)\right].\label{A}
\end{align}
\end{theorem}
\dproof \newline
Choose $\varphi \in C^{2}\left( \mathbb{R}\right) $ with bounded derivatives. Then by the It\^{o}
formula (Theorem \ref{Ito} ), we have%
\[
\mathbb{E}\left[ \varphi \left( X_{t+h}\right) -\varphi \left( X_{t}\right) %
\right] =\mathbb{E}\left[ \int_{t}^{t+h}A\varphi \left( X_{s}\right) ds%
\right] , 
\]%
where
\begin{align}
&A\varphi \left( X_s \right) =\alpha \left( s,X_s,\mu_s\right) \varphi ^{\prime
}\left( X_{s}\right) +\beta \left( s,\mu_s\right)M^2(\beta \chi_{[0,s]})(s) \varphi
^{^{\prime \prime }}\left( X_{s}\right) 
\end{align}
In particular, choosing, with $i=\sqrt{-1}$,
\begin{align}
\varphi \left(x \right) =\varphi _{y}\left( x\right) =e^{-iyx};\quad y,x\in 
\mathbb{R},
\end{align}%
we get%
\begin{align} 
A\varphi _{y}\left( X_{s}\right) 
 =\Big( -iy\alpha \left( s,X_{s},\mu_s\right) -%
y^2\beta \left( s,\mu_s\right)M^2(\beta \chi_{[0,s]})(s)  \Big) 
 e^{-iyX_{s}}.\label{3.4}
\end{align}
In general we have 
\[
\mathbb{E}\left[ g\left( X_{s}\right) e^{-iyX_{s}}\right] =\int_{\mathbb{R}%
}g\left( x\right) e^{-iyx}\mu _{s}\left( dx\right) =F\left[ g\left( \cdot\right)
\mu _{s}(\cdot)\right] \left( y\right).
\]
Recall that if $w \in \mathcal{S}'$ with $\left( \frac{d}{dx}%
\right) ^{n}w\left( t,x\right) =:D^{n}w\left( t,x\right) ,$ then we have, in the sense of distributions,
\[
F\left[ D^n w\left( t,\cdot  \right) \right] \left( y\right)
=\left( iy\right) ^{n}F\left[ w\left( t,\cdot\right) \right] \left( y\right). 
\]
Therefore
\begin{align}\label{D-alpha}
&-iyF[\alpha(s. \cdot)\mu_s](y)= F[-D(\alpha(s,\cdot)\mu_s)](y)\\\label{D-beta}
&-\tfrac{1}{2} y^2F[\beta^2(s,\cdot)\mu_s](y)=F[\tfrac{1}{2}D^2(\beta^2(s,\cdot) \mu_s)](y).
\end{align}
Applying this to \eqref{3.4} we get 
\begin{align}
\mathbb{E}[A\varphi _{y}( X_{s}) ]  &=-iyF[
\alpha ( s,\cdot,\mu_s) \mu _{s}(\cdot)] ( y) -
y^{2}F[ \beta ( s,\mu_s)M^2(\beta \chi_{[0,s]})(s)\mu_s(\cdot)] (y).  \label{L1}
\end{align}
Hence
\begin{align}\label{1.8}
&\widehat{\mu }_{t+h}\left( y\right) -\widehat{\mu }_{t}\left( y\right) 
=\int_{\mathbb{R}}e^{-iyx}\mu _{t+h}\left( dx\right) -\int_{\mathbb{R}%
}e^{-iyx}\mu _{t}\left( dx\right) \nonumber \\
&=\mathbb{E}[ \varphi _{y}\left( X_{t+h}\right) -\varphi _{y}\left(
X_{t}\right) ] =\mathbb{E}[ \int_t^{t+h}A\varphi _{y}( X_{s}) ds ] =\int_{t}^{t+h}L_{s}\left( y\right) ds,
\end{align}
where
\begin{align}\label{L1}
L_s(y)&=\mathbb{E}\left[ A\varphi _{y}\left( X_{s}\right) \right]\\ &=-iyF\left[
\alpha \left( s,\cdot,\mu_s\right) \mu _{s}(\cdot)\right] \left( y\right) -
y^{2} \beta ( s,\mu_s)M^2(\beta \chi_{[0,s]})(s)F\left[\mu_s(\cdot)\right] \left( y\right)\nonumber \\
\end{align}
Since $s \mapsto L_s$ is continuous, we conclude from \eqref{1.8} that
\begin{align} \label{3.8}
\tfrac{d}{dt}\widehat{\mu}_t(y)= L_t(y);\quad y \in \mathbb{R}.
\end{align}
Combining \eqref{D-alpha},\eqref{D-beta},  we get
\begin{align}\label{3.13}
L_s(y)= F[A^{*}\mu_s](y),
\end{align}
where $A^{*}$ is the differential operator
\begin{align}
A^{*}\mu_s =-D[\alpha(s,\cdot,\mu_s)\mu_s(\cdot)]+\beta ( s,\mu_s)M^2(\beta \chi_{[0,s]})(s)D^2[\mu_s(\cdot)].
\end{align}
Note that $A^{*}\mu_s$ exists in $\mathcal{S}^{\prime}$.\\
Combining \eqref{3.8} with \eqref{3.13}, we get
\begin{align}
\tfrac{d}{dt} F[\mu_t](y)= F[A^{*}\mu_t](y).
\end{align}
Hence
\begin{align}
F[\mu_t](y)= F[\mu_0](y)+ \int_0^t F[A^{*}\mu_s](y) ds=F[\mu_0 + \int_0^t (A^{*}\mu_s)ds](y).
\end{align}
Since the Fourier transform of a distribution determines the distribution uniquely, we deduce that
\begin{align}
\mu_t= \mu_0 + \int_0^t A^{*}\mu_s ds \text{ in } \mathcal{S^{\prime}}.
\end{align}
Therefore 
\begin{align}
\tfrac{d}{dt} \mu_t= A^{*}\mu_t \text{ in } \mathcal{S^{\prime}},
\end{align}
as claimed.
\fproof 

\begin{corollary}[Weak Fokker--Planck equation for the density]
Assume in addition that for each $t\in[0,T]$ the law $\mu_t$ of $X_t$ admits a density
$m_t\in L^1(\R)$, i.e.\ $\mu_t(dx)=m_t(x)\,dx$.
Then the map $t\mapsto m_t$ is absolutely continuous as an $\mathbb{M}_0$-valued
(and hence $\mathcal S'$-valued) curve, and for all $\phi=\phi(t,x)\in\mathcal S(\R\times \R)$,
\begin{align}
    &\langle \frac{d}{dt} \mu_t, \phi(t,x)\rangle = \langle A^{*}\mu_t,\phi(t,x)\rangle \  \text{ i.e. }\\
    &-\int_{\R} \frac{d}{dt}\phi(t,x) \mu_t(dx)= \int_{\R} A\phi(t,x) \mu_t(dx).
\end{align}

In particular, $m$ is a weak solution of the Fokker--Planck equation
\[
\partial_t m_t(x)
= -\partial_x\big(\alpha(t,x,\mu_t)m_t(x)\big)
  + \beta(t,\mu_t)M^2(\beta\chi_{[0,t]})(t)\,\partial_x^2 m_t(x)
\]
on $(0,T)\times\R$.
\end{corollary}
\begin{remark}
We do not assume any $C^{1,2}$ regularity of $(t,x)\mapsto m_t(x)$.
The time derivative $\partial_t m_t$ in the above equation is understood
in the distributional sense via \eqref{eq:weak-FP-density}.
If, in addition, $m\in C^{1,2}([0,T]\times\R)$, then the weak formulation
implies that $m$ is a classical solution of the Fokker--Planck equation.
\end{remark}

\section{A Feynman--Kac Representation for Fractional McKean--Vlasov Diffusions}
One of the most fundamental links between stochastic processes and partial
differential equations is provided by the \emph{backward Kolmogorov equation},
also known as the \emph{Feynman--Kac formula}. 
In classical diffusion theory, this identity expresses solutions of a
parabolic PDE in terms of expectations of a stochastic process.  
In the present setting, we show that an analogous result holds for the
fractional McKean--Vlasov diffusion \eqref{MK} driven by fBm with Hurst parameter $H>1/2$.
More precisely, the Fokker--Planck equation established in the previous
section gives the forward-in-time evolution of the law $\mu_t$ of $X_t$.  
The goal of this section is to derive the \emph{dual}, backward equation for
functionals of the form
\[
v(t,x) := \E^x\!\left[\,\varphi(X_t)\,\right],
\]
which represent the expected value at time $t$ of a terminal payoff
$\varphi$ when the process is started from $x$ at time $0$.  
We show that $v$ solves the PDE
\[
\partial_t v = A v,
\]
where $A$ is the infinitesimal generator associated with the fractional
McKean--Vlasov dynamics.  
This is the fractional analogue of the classical Feynman--Kac representation.
\vspace{0.25cm}
\begin{theorem}[Feynman--Kac formula for the fractional McKean--Vlasov SDE]
\label{thm:FK}
Let $X(t)$ be the solution of the McKean--Vlasov fBm equation \eqref{MK}, 
and assume that the law of $X_t$ admits a density,
\[
\mu_t(dx) = m_t(x)\,dx,
\]
for each $t\in[0,T]$.
Let $\varphi:\R\to\R$ be a smooth function with bounded support and define
\[
v(t,x) := \E^x[\varphi(X_t)] , \qquad t\ge 0,\ x\in\R .
\]
Let $A, A^*$ denote the generators appearing in the Fokker--Planck equation
(see \eqref{A},\eqref{A*}), so that
\[
\langle A\phi ,\, \mu\rangle = \langle \phi,\, A^*\mu\rangle; \  \phi \in \mathcal{S}(\R).
\]
Then $v$ satisfies the backward Kolmogorov equation
\begin{equation}\label{eq:BK}
\frac{\partial}{\partial t} v(t,x)
  = A v(t,\cdot)(x),  
\qquad t>0,\ x\in\R,
\end{equation}
with initial condition $v(0,x)=\varphi(x)$.
\end{theorem}
\dproof
Fix $x\in\R$ and let $\mu^x_t$ denote the law of $X_t$ under $\P^x$.  
Since $\mu^x_t(dx)=m^x_t(x)\,dx$, we may write
\[
v(t,x) 
 = \E^x[\varphi(X_t)]
 = \int_{\R} \varphi(r)\,m_t^x(r)\,dr .
\]
Using Theorem~\ref{FP1} (Fokker--Planck equation), we have, in the sense of
distributions,
\[
\partial_t m_t^x = A^* m_t^x ,
\]
and therefore
\[
\frac{\partial}{\partial t}v(t,x)
 = \int_{\R} \varphi(r)\,\partial_t m_t^x(r)\,dr
 = \int_{\R} \varphi(r)\,A^* m_t^x(r)\,dr .
\]
By definition of the adjoint,
\[
 \int_{\R} \varphi(r)\,A^* m_t^x(r)\,dr
 = \int_{\R} (A\varphi)(r)\,m_t^x(r)\,dr
 = \E^x\!\left[(A\varphi)(X_t)\right].
\]
Finally, as $A$ acts only on the spatial variable,
\[
\frac{\partial}{\partial t} v(t,x)
 = \E^x\!\left[(A\varphi)(X_t)\right]
 = A\Big(\E^x[\varphi(X_t)]\Big)(x)
 = A v(t,x),
\]
which proves \eqref{eq:BK}.  
The initial condition is immediate from $X_0=x$.
\fproof
\section{Examples}
In this section, we apply our techniques to solve two examples, demonstrating their effectiveness and utility.
\subsection{The law of fractional Brownian motion}
To illustrate our general results, we recover the marginal law of
fBm from the associated Fokker--Planck equation.
Consider the fractional SDE
\[
dX(t)=dB^{(H)}(t), \qquad X(0)=0.
\]
This corresponds to choosing $Z=0$, $\alpha=0$ and $\beta=1$ in the general
framework.
The associated Fokker--Planck equation \eqref{FP1} takes the form
\begin{equation}\label{fp4}
\partial_t m_t(x)
= \lambda(t)\,\partial_x^2 m_t(x),
\qquad t>0,\ x\in\R,
\end{equation}
where
\[
\lambda(t)=M^2(\chi_{[0,t]})(t).
\]
By \eqref{M2}, this can be written as
\[
\lambda(t)
= 2H C_H t^{2H-1}.
\]
Equation \eqref{fp4} is a heat equation with time dependent diffusivity $\lambda(t)$.
Its solution is given by
\begin{equation}\label{m1}
m_t(x)=p(\Lambda(t),x), \text{ with } \Lambda(t)=\int_0^t \lambda (s) ds= C_H t^{2H},
\end{equation}
where 
\[
p(u,x)=\frac{1}{\sqrt{4\pi u}}
\exp\!\left(-\frac{x^2}{4u}\right).
\]
Hence
\begin{align}
    m_t(x)=\frac{t^{H}}{(4 \pi C_H)^{\frac{1}{2}}} \exp(-\frac{x^2}{4C_H t^{2H}}); \ t > 0.\label{m_t}
\end{align}
\begin{remark}
    We can also obtain this directly as follows:\\
Since fBm is a
centered Gaussian process with variance
\[
\E[B^{(H)}(t)^2]=2C_H t^{2H},
\]
(see \eqref{var}), its density is
\begin{equation}\label{m2}
m_t(x)
=\frac{1}{\sqrt{4\pi C_H t^{2H}}}
\exp\!\left(-\frac{x^2}{4C_H t^{2H}}\right); \ t>0.
\end{equation}
which coincides with \eqref{m_t}.
\end{remark}

 \begin{figure}[!htb]
\centering
\begin{subfigure}{.5\textwidth}
  \centering
  \includegraphics[width=1\linewidth]{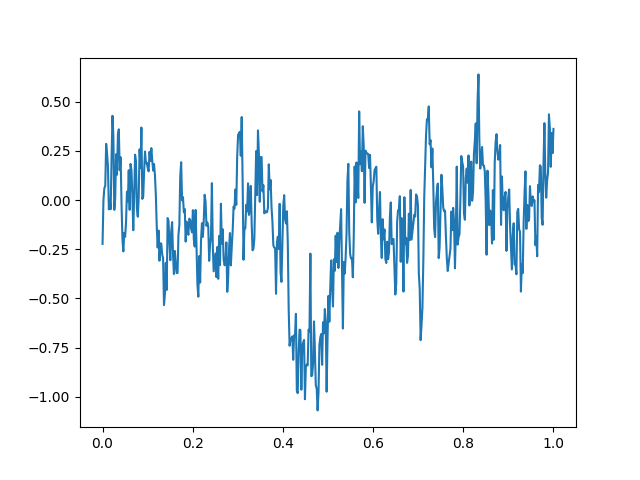}
 \caption{\textbf{Realisation for $H=0.3$}}
\end{subfigure}%
\begin{subfigure}{.5\textwidth}
  \centering
  \includegraphics[width=1\linewidth]{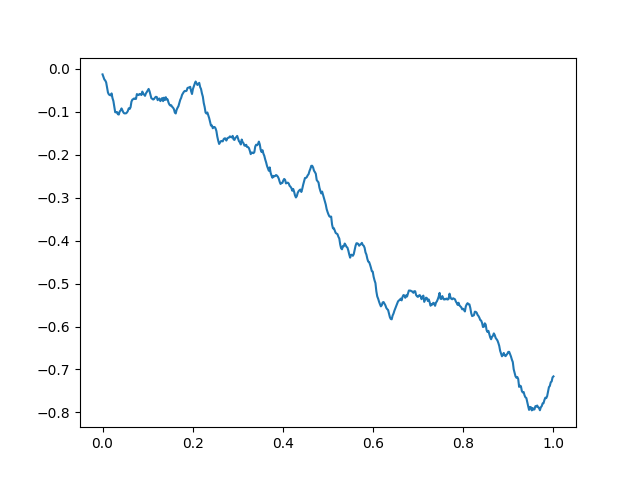}
 \caption{\textbf{Realisation for $ H=0.8$}}
\end{subfigure}
\end{figure}
 \begin{figure}[!htb]
\centering
\begin{subfigure}{.5\textwidth}
  \centering
  \includegraphics[width=1\linewidth]{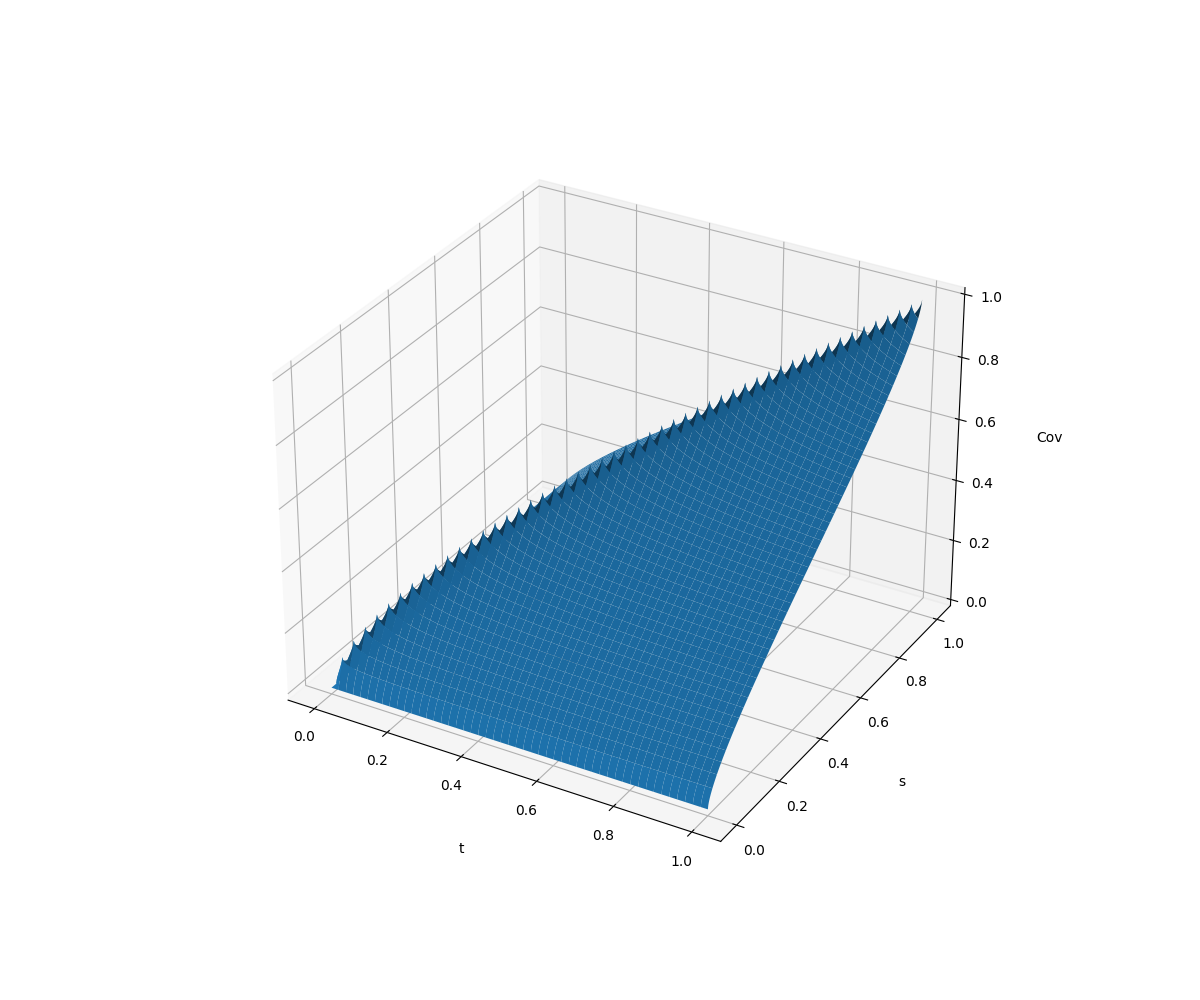}
 \caption{$Cov(B^{(H)}(t),B^{(H)}(s))_{H=0.3}$}
\end{subfigure}%
\begin{subfigure}{.5\textwidth}
  \centering
  \includegraphics[width=1\linewidth]{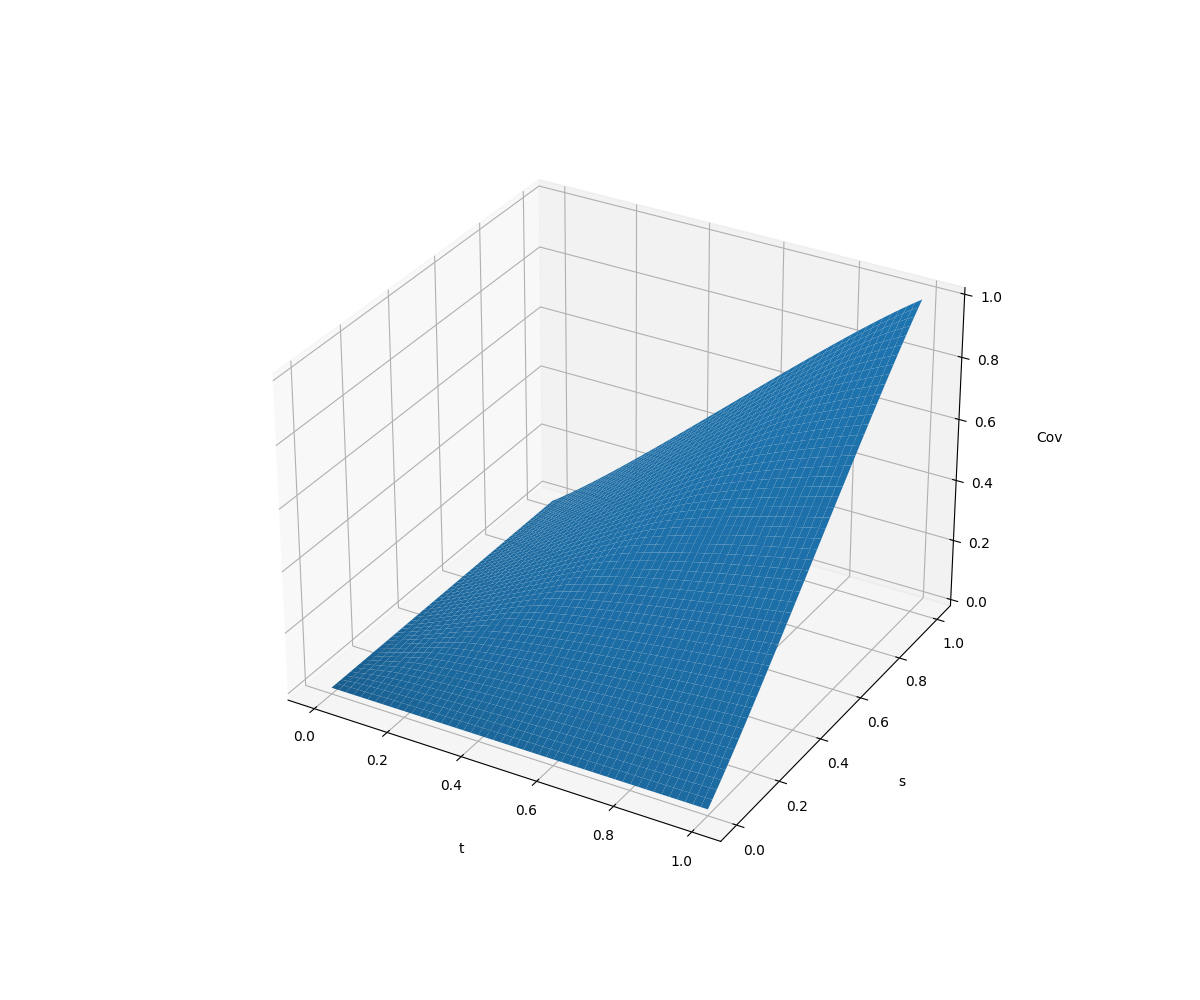}
 \caption{$Cov(B^{(H)}(t),B^{(H)}(s))_{H=0.8}$}
\end{subfigure}
\end{figure}
\subsection{The law of a linear McKean--Vlasov fractional SDE}
We consider the following linear McKean--Vlasov equation driven by fBm:
\begin{align}\label{eqX}
dX(t)
= \alpha_0 \E[X(t)]\,dt + \beta_0 \E[X(t)]\,dB^{(H)}(t),
\qquad X(0)=Z,
\end{align}
where $\alpha_0,\beta_0\in\R$ with $\beta_0\neq0$ and $Z\in L^2(\Omega)$ is
$\mathcal F_0$-measurable.
Here the coefficients are given by
\[
\alpha(t,x,\mu)=\alpha_0\int_{\R} r\,\mu(dr),
\qquad
\beta(t,\mu)=\beta_0\int_{\R} r\,\mu(dr).
\]
\paragraph{Fokker--Planck equation.}
Let $\mu_t(dx)=m_t(x)\,dx$ denote the law of $X(t)$.  
The associated Fokker--Planck equation reads
\[
\partial_t m_t(x)
= -\alpha_0\Big(\int_{\R} r\,m_t(dr)\Big)\partial_x m_t(x)
 + \beta_0^2\Big(\int_{\R} r\,m_t(dr)\Big)^2
   M^2(\chi_{[0,t]})(t)\,\partial_x^2 m_t(x).
\]
Taking expectations in \eqref{eqX} yields
\[
\frac{d}{dt}\E[X(t)] = \alpha_0 \E[X(t)],
\qquad
\E[X(0)]=\E[Z],
\]
and therefore
\[
\E[X(t)] = \E[Z]\,e^{\alpha_0 t}.
\]
Substituting this into the Fokker--Planck equation gives
\begin{align}\label{fp}
\partial_t m_t(x)
= -\alpha_0 \E[Z]e^{\alpha_0 t}\,\partial_x m_t(x)
 + \beta_0^2 \E[Z]^2 e^{2\alpha_0 t}
   M^2(\chi_{[0,t]})(t)\,\partial_x^2 m_t(x).
\end{align}
\paragraph{Explicit law and Feynman--Kac representation.}
Equation \eqref{eqX} is equivalent to
\[
dX(t)
= \alpha_0 x_0 e^{\alpha_0 t}\,dt
 + \beta_0 x_0 e^{\alpha_0 t}\,dB^{(H)}(t),
\qquad x_0:=\E[Z].
\]
Hence
\[
X(t)
= Z
 + x_0\frac{e^{\alpha_0 t}-1}{\alpha_0}
 + \beta_0 x_0 \int_0^t e^{\alpha_0 s}\,dB^{(H)}(s).
\]
Let $\varphi\in C_c^\infty(\R)$ and define
\[
v(t,x)=\E^x[\varphi(X(t))].
\]
By the fractional Feynman--Kac theorem proved above, $v$ solves the backward
equation
\begin{align}\label{eq:BK-geometric}
\partial_t v(t,x)
&= \alpha_0 x_0 e^{\alpha_0 t}\,\partial_x v(t,x)
 + \beta_0^2 x_0^2 e^{2\alpha_0 t}
   M^2(\chi_{[0,t]})(t)\,\partial_x^2 v(t,x),\\
v(0,x)&=\varphi(x).\nonumber
\end{align}
Moreover, the solution admits the explicit stochastic representation
\[
v(t,x)
= \E\!\left[
\varphi\!\left(
x + x_0(e^{\alpha_0 t}-1)
+ \beta_0 x_0\int_0^t e^{\alpha_0 s}\,dB^{(H)}(s)
\right)
\right].
\]
In particular, $X(t)$ is Gaussian with explicitly known mean and variance,
and its law is fully characterized.  Equation \eqref{fp} is therefore a
time inhomogeneous fractional heat equation with explicit coefficients.

\end{document}